\documentclass[12pt]{article}

\newcommand{\mathR}{{\mathbb R}}
\newcommand{\mathZ}{{\mathbb Z}}

\usepackage{graphics}
\usepackage{amsmath,amssymb}
\usepackage{graphicx,latexcad}
\usepackage[dvips]{xy}
\xyoption{all}

\begin{document}
\baselineskip=6.0mm

\newcommand{\ignore}[1]{}{}

\renewcommand{\theequation}{\arabic{section}.\arabic{equation}}

\newcommand{\lbl}{\label}


\newcommand{\eq}[1]{$(\ref{#1})$}

\newcommand{\al}{\alpha}                         
\newcommand{\bt}{\beta}                          
\newcommand{\ga}{\gamma}                         
\newcommand{\Ga}{\Gamma}                         
\newcommand{\de}{\delta}                         
\newcommand{\De}{\Delta}                         
\newcommand{\ep}{\epsilon}                       
\newcommand{\ve}{\varepsilon}                    
\newcommand{\la}{\lambda}                        
\newcommand{\La}{\Lambda}                        
\newcommand{\ro}{\rho}                           
\newcommand{\ta}{\tau}                           
\newcommand{\si}{\sigma}                         

\newcommand{\be}{\begin{equation}}               
\newcommand{\ee}{\end{equation}}                 
\newcommand{\bea}{\begin{eqnarray}}              
\newcommand{\eea}{\end{eqnarray}}                
\newcommand{\bean}{\begin{eqnarray*}}            
\newcommand{\eean}{\end{eqnarray*}}              
\newcommand{\ba}{\begin{array}}                  
\newcommand{\ea}{\end{array}}                    
\newcommand{\nn}{\nonumber}                      
\newcommand{\mb}{\mbox}                          

\newcommand{\ra}{\rightarrow}                    
\newcommand{\llra}{\longleftrightarrow}          

\newcommand{\stac}{\stackrel}                    
\newcommand{\noin}{\noindent}                    

\newcommand{\qed}{\nobreak\quad\vrule width6pt depth3pt height10pt}

\newcommand{\ngi}{n \ra \infty}

\pagestyle{myheadings} \markright{Rates of convergence of means of
Euclidean functionals}

\thispagestyle{plain}

\begin{center}
{\Large\bf Rates of convergence of means of Euclidean functionals}
\vskip 0.4in

Yooyoung Koo,$^{a,}$\footnote{Supported by the BK21 project of the
Department of Mathematics, Sungkyunkwan University.} and Sungchul
Lee,$^{b,}$\footnote{Supported by the BK21 project of the
Department of Mathematics, Yonsei University.}

\vskip 0.2in

$^a$Department of Mathematics,
Sungkyunkwan University, Suwon, Korea
440-746

$^b$Department of Mathematics, Yonsei University, Seoul, Korea
120-749

\vskip 0.3in
\end{center}
\vskip 0.3in

\begin{abstract}
Let $L$ be the Euclidean functional with $p$-th power-weighted
edges. Examples include the sum of the $p$-th power-weighted
lengths of the edges in minimal spanning trees, traveling salesman
tours, and minimal matchings. Motivated by the works of Steele,
Redmond and Yukich (1994, 1996) have shown that for $n$ i.i.d.\ sample points $\{X_1,\ldots,X_n\}$ from $[0,1]^d$,
$L(\{X_1,\ldots,X_n\})/n^{(d-p)/d}$ converges a.s.\ to a finite
constant. Here we bound the rate of convergence of
$EL(\{X_1,\ldots,X_n\})/n^{(d-p)/d}$.

\vspace{1.0cm}
\hrule
\vspace{0.1cm}
\noindent {\bf \it AMS $1991$ subject classification.}
                           Primary 60D05;
                           secondary 05C80, 90C27.
\vspace{0.1cm}
\hrule
\end{abstract}

\section{Introduction.}
\setcounter{equation}{0}

Let $\{X_1,\ldots,X_n\}$ be $n$ i.i.d.\ sample points from
$\mathR^d$, $d \geq 2$, and let $0<p<\infty$.  A traveling
salesman problem (TSP) is to find a permutation $\pi$ on
$\{1,\ldots,n\}$ such that \bea
& & \sum_{j=1}^{n} |X_{\pi(j+1)}-X_{\pi(j)}|^p \nn \\
&=& \min  \left\{ \sum_{j=1}^{n} |X_{\pi'(j+1)}-X_{\pi'(j)}|^p
 : \pi' \mbox{ a permutation on $\{1,\ldots,n\}$}  \right\}, \nn
\eea where $|X_i-X_j|$ is the Euclidean distance between $X_i$ and
$X_j$ and where $\pi(n+1):=\pi(1)$ and $\pi'(n+1):=\pi'(1)$. Let $
L_{TSP}\left(\{X_1,\ldots,X_n\},p\right)$ be the sum of the $p$-th
power-weighted lengths of the edges in a minimal tour $\pi$. In
the case $\{X_1,\ldots,X_n\}=\emptyset$ define
$L_{TSP}(\emptyset,p)=0.$ Beardwood, Halton, and Hammersley (1959)
showed that there exists a strictly positive but finite constant
$\alpha(L_{TSP},d,1)$ such that for i.i.d.\ sample points $\{
X_i : i \geq 1 \}$ with common distribution $\mu$, which has a
compact support in $\mathR^d$, $d \geq 2$, as $n \ra \infty$ \be
\frac{L_{TSP}(\{X_1,\ldots,X_n\},1)}
     {n^{(d-1)/d}}
\rightarrow \alpha(L_{TSP},d,1) \int f^{(d-1)/d}(x) dx \mbox{
a.s.} \lbl{1.1} \ee where $f$ is the density function of the
absolutely continuous part of $\mu$.

The asymptotic behavior \eq{1.1} of the TSP functional is not an
isolated one. A minimal matching (MM) on $\{X_1,\ldots,X_n\}$ is a
permutation $\pi$ on $\{1,\ldots,n\}$ such that \bea
&& \sum_{j=1}^{[n/2]} |X_{\pi(2j)}-X_{\pi(2j-1)}|^p \nn\\
&=& \min  \left\{ \sum_{j=1}^{[n/2]}
|X_{\pi'(2j)}-X_{\pi'(2j-1)}|^p
 : \pi' \mbox{ a permutation on $\{1,\ldots,n\}$}  \right\}, \nn
\eea where $[n/2]$ is the largest integer $l$ with $l \leq n/2$.
Let $ L_{MM}\left(\{X_1,\ldots,X_n\},p\right) $ be the sum of the
$p$-th power-weighted lengths of the edges in a minimal matching
$\pi$. In the case $\{X_1,\ldots,X_n\}=\emptyset$ define
$L_{MM}(\emptyset,p)=0.$

A minimal spanning tree (MST) on $\{X_1,\ldots,X_n\}$ is a
spanning tree $T$ on the given point set $\{X_1,\ldots,X_n\}$ such
that \bea &&
\sum_{(X_i,X_j) \in T} |X_i-X_j|^p \nn\\
&=& \min  \left\{ \sum_{(X_i,X_j) \in T'} |X_i-X_j|^p : T' \mbox{
a spanning tree on $\{X_1,\ldots,X_n\}$}  \right\}.\nn \eea Let $
L_{MST}\left(\{X_1,\ldots,X_n\},p\right) $ be the sum of the
$p$-th power-weighted lengths of the edges in a minimal spanning
tree $T$. In the case $\{X_1,\ldots,X_n\}=\emptyset$ define
$L_{MST}(\emptyset,p)=0.$

A Steiner minimal spanning tree (SMST) on $\{X_1,\ldots,X_n\}$ is
a spanning tree $T$ on a point set containing $\{X_1,\ldots,X_n\}$
(we call such $T$  a Steiner spanning tree  on
$\{X_1,\ldots,X_n\}$) such that \bea &&
\sum_{(X_i,X_j)\in T} |X_i-X_j|^p \nn\\
&=& \min  \left\{ \sum_{(X_i,X_j) \in T'} |X_i-X_j|^p : T' \mbox{
a Steiner spanning tree on $\{X_1,\ldots,X_n\}$} \right\}.\nn \eea
Let $ L_{SMST}\left(\{X_1,\ldots,X_n\},p\right) $ be the sum of
the $p$-th power-weighted lengths of the edges in a minimal
Steiner spanning tree $T$. In the case
$\{X_1,\ldots,X_n\}=\emptyset$ define $L_{SMST}(\emptyset,p)=0.$
Note that for $p > 1$, $L_{SMST}=0$. So, whenever we talk about
$L_{SMST}$, we always consider the case $0<p\le 1$.

A rectilinear Steiner minimal spanning tree (RSMST) on
$\{X_1,\ldots,X_n\}$ is a Steiner spanning tree $T$ on
$\{X_1,\ldots,X_n\}$ in which all the edges are rectilinear (we
call such $T$  a rectilinear Steiner spanning tree  on
$\{X_1,\ldots,X_n\}$) such that \bea &&
\sum_{(X_i,X_j)\in T} |X_i-X_j|^p \nn\\
&=& \min  \left\{ \sum_{(X_i,X_j) \in T'} |X_i-X_j|^p : T' \mbox{
a rectilinear spanning tree on $\{X_1,\ldots,X_n\}$} \right\}.\nn
\eea Let $ L_{RSMST}\left(\{X_1,\ldots,X_n\},p\right) $ be the sum
of the $p$-th power-weighted lengths of the edges in a minimal
rectilinear Steiner spanning tree $T$. In the case
$\{X_1,\ldots,X_n\}=\emptyset$ define $L_{RSMST}(\emptyset,p)=0.$
Note again that for $p > 1$, $L_{RSMST}=0$. So, whenever we talk
about $L_{RSMST}$, we always consider the case $0<p\le 1$.

In a series of papers Steele (1981A, 1981B, 1988, 1990) showed
that the asymptotic behavior \eq{1.1} appears for various
functionals including some of the above five. Redmond and Yukich
(1994, 1996) and Lee (1999) further developed the general conditions
providing the asymptotic behavior \eq{1.1}. One may consult
Section 1.2 of Yukich (1998) for a brief history of this field.

\ignore{ Then, Papadimitriou (1978) and Steele (1988) showed that
$L_{MM}$ and $L_{MST}$ have the same asymptotics as that of
$L_{TSP}$ where $\alpha(L_{TSP},d,1)$ is replaced by
$\alpha(L_{MM},d,1)$ and $\alpha(L_{MST},d,1)$, respectively.
More generally Redmond and Yukich (1994) found a sufficient
condition for the functional $L$ which guarantees the asymptotics
\eq{1.1} for $L$.  Redmond and Yukich (1996) developed their
theory further and showed that for $L=L_{MM}$, $L_{MST}$,
$L_{TSP}$, and for $1 \leq p < d$, there exists a strictly
positive but finite constant $\alpha(L,d,p)$ such that for i.i.d.\
sample points $\{ X_i : i \geq 1 \}$ with common distribution
$\mu$, which has a compact support in $\mathR^d$, $d \geq 2$, as
$n \ra \infty$ \be \frac{L(\{X_1,\ldots,X_n\},p)}
     {n^{(d-p)/d}}
\rightarrow \alpha(L,d,p) \int f^{(d-p)/d}(x) dx \mbox{ c.c.}
\lbl{1.2} \ee where $f$ is the density function of the absolutely
continuous part of $\mu$.  Here, c.c.\ means the complete
convergence which is a stronger notion than the a.s.\ convergence.
One can find good surveys on this subject from Steele (1997) and
Yukich (1998).

In this thesis we prove in fact, for MM, TSP, MST, SMST and RSMST,
the rate of convergence problem. Redmond and Yukich (1994, 1996)
use an Euclidean functional pair $(L_{MST},L^*_{MST})$ to get an
upper bound for $|EL_{MST}({\cal
U}_n,[0,1]^d,p)/n^{(d-p)/d}-\alpha(L_{MST},d,p)|$, $1 \le p<d$.
However, for $0<p<1$ the Euclidean functional pair argument seems
not working. And if $L$ is the minimal matching or TSP functional,
then it is unclear whether the rate of convergence is satisfied.}

Our results are stated below.  But, first we would like to spell
out the restrictions on the Euclidean functional $L$. We call
$L({\cal A},B,p)$, ${\cal A}$ a finite subset of a box
$B=\prod_{i=1}^d[x_i,x_i+s]$, $x =(x_1,\ldots,x_d)\in \mathR^d$,
$0<s<\infty$, $d \geq 2$, $0<p<\infty$, {\it a subadditive
Euclidean functional (or a weak Euclidean functional) } of power
$p$ if the following four conditions are met: \be
L(\emptyset,B,p)=0, \lbl{1.2} \ee for $y \in \mathR^d$ and
$0<t<\infty$ \be L(y+t{\cal A},y+tB,p)=t^pL({\cal A},B,p),
\lbl{1.3} \ee \be \left|L({\cal A},B,p)-L({\cal B},B,p)\right|
\leq C |{\cal A} \triangle {\cal B}|^{(d-p)/d}s^p, \lbl{1.4} \ee
and for a partition $\{Q_i, 1 \leq i \leq m^d\}$ of $[0,1]^d$ into
$m^d$ subboxes of edge length $m^{-1}$ \be L({\cal A},[0,1]^d,p)
\leq \sum_{j=1}^{m^d} L({\cal A} \cap Q_i,Q_i,p) +Cm^{d-p}.
\lbl{1.5} \ee

By \eq{1.2} and \eq{1.4} with $s=1$ and ${\cal B}=\emptyset$, we
have for a finite subset ${\cal A}$ of the unit box $[0,1]^d$ and
for $d \geq 2$, $0<p<d$, \bea |L({\cal A},[0,1]^d,p) \leq C |{\cal
A}|^{(d-p)/d}. \lbl{a} \eea More strongly, for all the Euclidean
functionals $L$ of our interest in this paper (look at Theorem 2
below for the full list of such $L$) based on the space filling
curve heuristic, as shown by Steele (1997), there is an extension
of the above bound \eq{a} which covers both the $d=1$ case and the
$p \geq d$ case; for a finite subset ${\cal A}$ of the unit box
$[0,1]^d$ and for $d \geq 1$, $p>0$, \bea |L({\cal A},[0,1]^d,p)|
\leq C(|{\cal A}|^{(d-p)/d} \vee 1). \lbl{*} \eea Note that \eq{a}
follows from the assumptions of the Euclidean functional whereas
\eq{*} follows from the specific feature of the Euclidean
functional of our interest in this paper.

Let ${\cal U}_n$ be $n$ i.i.d.\ uniform
points in $[0,1]^d$.\\

\noin {\bf Theorem A.} [Redmond and Yukich (1994, 1996),
Lee(1999)] Let $L$ be a weak Euclidean functional of power
$0<p<d$. Then there exists a finite constant $\al:=\alpha(L,d,p)$
such that as $n \ra \infty$ \be \frac{L({\cal
U}_n,[0,1]^d,p)}{n^{(d-p)/d}} \rightarrow \alpha \mbox{ c.c.\ and
in $L^1$,} \lbl{1.6} \ee where $Y_n \rightarrow \alpha$ c.c.\ (complete convergence) means
that for any $\ve > 0$, $\sum_{n=1}^{\infty} P(|Y_n -
\alpha| > \ve) < \infty$.\\

For a typical weak Euclidean functional $L$, the limit $\alpha$ in
\eq{1.6} is strictly positive:
In most situations of interest, the limit $\al$ is just the subadditive constant
and therefore must be strictly positive.

We call $L^*({\cal A},B,p)$ {\it a superadditive Euclidean
functional} of power $p$ if $-L^*({\cal A},B,p)$ is a subadditive
Euclidean functional of power $p$. We call $L$ {\it a Euclidean
functional} of power $p$ if $L$ is a subadditive Euclidean
functional of power $p$, if $L^*$ is a superadditive Euclidean
functional of power $p$, and if \be L^*({\cal A},[0,1]^d,p) \leq
L({\cal A},[0,1]^d,p), \lbl{1.7} \ee and for the $n$ uniform
points ${\cal U}_n$ in $[0,1]^d$ \be \left|EL({\cal
U}_n,[0,1]^d,p)-EL^*({\cal U}_n,[0,1]^d,p)\right|=o(n^{(d-p)/d}).
\lbl{1.8} \ee

\noin {\bf Theorem B.} [Redmond and Yukich (1994, 1996),
Lee(1999)] Let $L$ be a Euclidean functional of power $0<p<d$.
Then for i.i.d.\ sample points $\{ X_i : i \geq 1 \}$ with
common distribution $\mu$, which has a compact support in
$[0,1]^d$, as $n \ra \infty$ \be
\frac{L(\{X_1,\ldots,X_n\},[0,1]^d,p)}{n^{(d-p)/d}}\rightarrow
\alpha \int_{[0,1]^d} f^{(d-p)/d}(x) dx \mbox{ c.c.\ and in $L^1$},
\lbl{1.9} \ee where $\al:=\al(L,d,p)$ is a finite constant given
by \eq{1.6} and
where $f$ is the density function of the absolutely continuous part of $\mu$.\\

If there is no confusion, to save the heavy notations
 from now on in the case $B=[0,1]^d$
we use the notation $L({\cal A})$ and $L^*({\cal A})$ instead of
$L({\cal A},[0,1]^d,p)$ and $L^*({\cal A},[0,1]^d,p)$,
respectively.

We call $L$ {\it a strong Euclidean functional} of power $0<p<d$
if $L$ is a Euclidean functional of power $p$ and if for  the
homogeneous Poisson point process ${\cal P}_n$ of intensity $n$ on
$[0,1]^d$ \be \left|EL({\cal U}_n)-EL^*({\cal U}_n)\right|\le
\left\{ \ba{ll} C n^{(d-1-p)/d}&\mb{ for $0<p<d-1$} \\
C(\log n \vee 1) &\mb{ for $p=d-1\neq 1$}\\
C &\mb{ for $d-1<p<d$} \\
C &\mb{ for $p=d-1=1$}, \ea \right. \lbl{1.10} \ee \be
\left|EL({\cal P}_n)-EL^*({\cal P}_n)\right| \le
\left\{ \ba{ll} C n^{(d-1-p)/d}&\mb{ for $0<p<d-1$} \\
C(\log n \vee 1) &\mb{ for $p=d-1\neq 1 $}\\
C &\mb{ for $d-1<p<d$} \\
C &\mb{ for $p=d-1=1$}, \ea \right. \lbl{1.11} \ee \bea
& & \left|EL({\cal U}_{n \pm k})-EL({\cal U}_n)\right| \nn \\
&\le& \left\{\ba{lll} C (k^{(d-1-p)/(d-1)}\vee 1)&\mb{ for $1 \le k\le n^{(d-1)/d}$}\\
C kn^{-p/d}&\mb{ for  $n^{(d-1)/d} \le k\le n/2,0 < p < d-1$}\\
C kn^{-(d-1)/d}&\mb{ for $n^{(d-1)/d} \le k\le n/2,d-1 \le p <
d$}. \ea \right. \lbl{1.12} \eea

We call $L$ {\it a very strong Euclidean functional} of power
$0<p<d$ if $L$ is a strong Euclidean functional of power $p$ and
if $L$ satisfies the add-one bound; \be \left| EL({\cal U}_{n+1})
- EL({\cal U}_{n})\right| \leq C n^{-p/d}. \lbl{1.13} \ee
Our first result is as follows.\\

\noin {\bf Theorem 1.} (i) For a strong Euclidean functional $L$
of power $0<p<d$, \be
\left. \ba{ll} - C n^{(d-1-p)/2(d-1)}\\
-C\\
-C\\
-C\ea \right\} \leq EL({\cal U}_n)-\alpha n^{(d-p)/d} \leq \left\{
\ba{ll}
C n^{(d-1-p)/d}&\mb{ for $0<p<d-1$}\\
C (\log n \vee 1)&\mb{ for $p=d-1\neq 1$}\\
C &\mb{ for $d-1<p<d$} \\
C &\mb{ for $p=d-1=1$,} \ea \right. \nn \ee where $\al:=\al(L,d,p)$
is a finite constant given by \eq{1.6}.

(ii) For a very strong Euclidean functional $L$ of power $0<p<d$,
\be - C n^{\frac{1}{2}-\frac{p}{d}} \leq
EL({\cal U}_n)-\alpha n^{(d-p)/d} \\
\leq  \left\{ \ba{ll}
C n^{(d-1-p)/d}&\mb{ for $0<p<d-1$}\\
C (\log n \vee 1) &\mb{ for $p=d-1\neq 1$}\\
C &\mb{ for $d-1<p<d$} \\
C &\mb{ for $p=d-1=1$,}\ea \right. \nn \ee
where $\al:=\al(L,d,p)$ is a finite constant given by \eq{1.6}.\\

To see the point of Theorem 1 we compare it with Theorem 5.2 of
Yukich (1998). In Theorem 5.2 of Yukich (1998), he showed that, if
$L$ satisfies the close in mean approximation, that is,
$
\left|EL({\cal U}_n)-EL^*({\cal U}_n)\right|=o(n^{(d-p)/d})
$,
and the add-one
bound, then \be |EL({\cal U}_{n}) - \alpha n^{(d-p)/d}| \leq
C(n^{(d-1-p)/d} \vee 1). \nn \ee The provided rate is quite
satisfactory. However, as we see in Theorem 2 below, many typical $L$ do not satisfy the add-one bound and hence we cannot apply
Theorem 5.2 of Yukich (1998) to those $L$. We wish to provide the
same rate without the add-one bound so that typical $L$ has the
provided rate of convergence. What we find in Theorem 1 is that
this task can be done by strengthening the close in mean
approximation to \eq{1.10} and by adding two extra conditions
\eq{1.11} and \eq{1.12}. We also find that for the case $p=d-1
\neq 1$ typical $L$ do not satisfy the other condition used in
Theorem 5.2 of Yukich (1998)
 and we properly fix it in \eq{1.10} and \eq{1.11}. We also take care the case
$0<p<1$ which was excluded in Theorem 5.2 of Yukich (1998).

Note that the bounds in Theorem 1 follow from the assumptions of
the Euclidean functionals. If we use the specific property of the
Euclidean functional, in some cases we can get much better bounds
and surprisingly we can even get a strictly positive lower bound
for large $n$; see Jaillet (1993), Rhee (1994), Yukich (1998),
Lee (2000) for these specific results.

\ignore{
\noin {\bf Theorem 1.2.} For a finite subset ${\cal A}$ of a box
$B=\prod_{i=1}^d[x_i,x_i+s]$, $x=(x_1,\ldots,x_d) \in \mathR^d$,
$0<s<\infty$, $d \geq 2$, $0<p<d$, define  $L({\cal
A},B,p):=L_{MM}({\cal A},p)$, $L_{MST}({\cal A},p),$
$L_{TSP}({\cal A},p)$. And define $L({\cal A},B,p):=L_{SMST}({\cal
A},p)$, $L_{RSMST}({\cal A},p)$ for $0<p\leq 1$. Then $L$ is a
subadditive Euclidean functional of power $p$ so that \eq{1.6} for
$L$ holds. Moreover there exists a superadditive Euclidean
functional $L^*$ of power $p$ such that $(L,L^*)$ is an Euclidean
functional pair of power $p$
so that \eq{1.9} for $L$ holds.\\
}

For a given subadditive 
Euclidean functional $L$, 
to find the rate of convergence of $EL({\cal U}_n)$ using Theorem
1 we have to construct a superadditive
Euclidean functional $L^*$ 
which puts us into the situations considered in  Theorem 1. There
may be several ways to construct such an $L^*$.  The successful
$L^*$ is the one which uses  the boundary freely.  For example
let's consider the traveling salesman problem.  Suppose that there
are $n$ cities ${\cal A}$ in $[0,1]^d$.  Then $L_{TSP}({\cal
A},[0,1]^d,1)$ is the total mileage to travel all the  $n$ cities.
Suppose however that there are ``free" ways along the boundary
$\partial [0,1]^d$ of $[0,1]^d$ in which the government pays the
gas.  In this case we may save some gas by traveling along the
boundary and $L_{TSP}^*({\cal A},[0,1]^d,1)$ is the total mileage.

In the case $1 \leq p < d$ following the idea of Redmond and
Yukich (1994, 1996) we construct superadditive Euclidean
functionals $L^*$ for the MM, MST, TSP, SMST, RSMST (of course for
the SMST, RSMST we consider the $p=1$ case only).  They are \be
L_{MM}^*({\cal A},B,p) := \min \left\{L_{MM}({\cal A} \cup {\cal
B},B,p): {\cal B} \mb{ a finite subset of $\partial B$}\right\},
\nn \ee \be L_{MST}^*({\cal A},B,p) := L_{MST}({\cal A},B,p)
\wedge \min \left\{\sum_j L_{MST}({\cal A}_j \cup
\{b_j\},B,p)\right\}, \nn \ee where the minimum is taken over the
partition $\{{\cal A}_j\}$ of ${\cal A}$ and $b_j \in \partial B$,
\be L_{TSP}^*({\cal A},B,p) := L_{TSP}({\cal A},B,p) \wedge \min
\left\{\sum_j \tilde L_{TSP}({\cal A}_j \cup
\{b_j,b'_j\},B,p)\right\}, \nn \ee where the minimum is taken over
the partition $\{{\cal A}_j\}$ of ${\cal A}$ and $b_j,b'_j \in
\partial B$ and where for a finite subset $\{X_1,\ldots,X_n\}$ of
$B$ with $|\{X_1,\ldots,X_n\} \cap
\partial B| \geq 2$ \be \tilde L_{TSP}(\{X_1,\ldots,X_n\},B,p)
:=\min \left\{
\sum_{j=1}^{n-1}\left|X_{\pi(j+1)}-X_{\pi(j)}\right|^p\right\},
\nn \ee where the minimum is taken over the permutation $\pi$ on
$\{1,\ldots,n\}$ such that $X_{\pi(1)}$, $X_{\pi(n)} \in
\partial B$.  Note that in the definition of $\tilde L_{TSP}$ the
sum is up to $n-1$ so that we travel free from $X_{\pi(n)}$ to
$X_{\pi(1)}$ along the boundary $\partial B$. Similarly, \be
L_{SMST}^*({\cal A},B,p) := L_{SMST}({\cal A},B,p) \wedge \min
\left\{\sum_j L_{SMST}({\cal A}_j \cup \{b_j\},B,p)\right\}, \nn
\ee where the minimum is taken over the partition $\{{\cal A}_j\}$
of ${\cal A}$ and $b_j \in
\partial B$,
\be L_{RSMST}^*({\cal A},B,p) := L_{RSMST}({\cal A},B,p) \wedge
\min \left\{\sum_j L_{RSMST}({\cal A}_j \cup \{b_j\},B,p)\right\},
\nn \ee where the minimum is taken over the partition $\{{\cal
A}_j\}$ of ${\cal A}$ and $b_j \in
\partial B$.
In the case ${\cal A}=\emptyset$ define $L_{MM}^*(\emptyset,B,p)=
L_{MST}^*(\emptyset,B,p)= L_{TSP}^*(\emptyset,B,p)=
L_{SMST}^*(\emptyset,B,p)= L_{RSMST}^*(\emptyset,B,p)=0. $

In the case $1 \leq p < d$ following the idea of Lee (1999) we
construct superadditive Euclidean functionals $L^*$ for the MM,
MST, TSP, SMST, and RSMST.  They are given in the following way. In
the above $L^*$ for $1 \le p <d$, there are some edges  $(X,Y)$
from a boundary point $X \in
\partial B$. In the case $0<p<1$, for any matching,
tree, or tour, we don't pay the full price  for the edge $(X,Y)$
from  the boundary point $X \in
\partial B$; for this edge we pay half of the
full price $\left|X-Y\right|^p$, i.e., $\left|X-Y\right|^{p}/2$. \\

\noin {\bf Theorem 2.} (i) For the Euclidean functional $L$ of MM,
MST, TSP with $0<p<d$ and  SMST, RSMST with $0<p\le 1$, $L$ is a
strong Euclidean functional (and hence  Theorems 1 (i), 3, 4 can
be applicable to this $L$).

(ii) For the Euclidean functional $L$ of MST with $0<p<d$ and
TSP, SMST, RSMST with $0<p\le 1$, $L$ is a very strong
Euclidean functional (and hence  Theorems 1 (ii), 3, 4 can be applicable to this $L$).\\

In Section 2, we develop a theory on the rate of convergence of
$EL$ for the uniform sample points, i.e., we prove Theorems 1 and
2. In Section 3, we continue to build a theory  on the rate of
convergence of $EL$ for the non-uniform sample points which was
started by Hero,  Costa,  and Ma (2003). \ignore{In Hero,  Costa,
and Ma (2003) they studied on the rate of convergence for $EL$ for
the non-uniform sample points  based on \eq{1.13} which is very
restrictive as you see at Theorem 2 (2) above. Our starting point
of this paper is to try to remove this condition in their argument
and our work in Section 2 is the result of this trial. In Section
3 we build  a theory  on the rate of convergence for $EL$ for the
non-uniform sample points  which does not depends on
 \eq{1.13} and furthermore
we extend the work of Hero,  Costa,  and Ma (2003) on $0<\bt \le
1$ to the general $\bt$ case.}

\section
{Rates of convergence; the uniform case} \setcounter{equation}{0}

In this section we prove Theorems 1 and 2. The main idea comes
from the symmetry argument of the patching in Alexander (1994).
Using this symmetry argument we get the nice moment estimate.
\ignore{ During the calculation on the patching we encounter some
moments which is nicely estimated by the symmetry argument.}
\\
\ignore{ Note that
$$
L({\cal P}_n,R,p) =_d n^{-1/d} L({\cal P},n^{1/d}R,p),
$$
where $R \subset \mathbb{R}^d$ is the set of sites. We will obtain
an analog of \eq{1.15} for the Poisson process \be \alpha
n^{(d-p)/d}-C \leq EL({\cal P}_n,[0,1]^d,p) \le \alpha
n^{(d-p)/d}+\left\{ \ba{ll}
Cn^{(d-1-p)/d}&\mb{ for $0<p<d-1$}\\
C\log n&\mb{ for $p=d-1$}\\
C&\mb{ for $d-1<p<d$}. \ea \right. \ee To obtain \eq{1.15}, we
first prove the above inequality. For fixed $m$ we divide the unit
cube into $m^d$ cubes of side $1/m$, and label these
$Q_1,\ldots,Q_{m^d}$. We also define
$$
\varphi (u):=EL({\cal P},[0,u]^d,p)$$ and $$\varphi^{*}
(u):=EL^*({\cal P},[0,u]^d,p).
$$
}

\noin {\bf Lemma 1.} Let $L$ be a strong Euclidean functional of
power $0<p<d$. Then,
$$
-C \le EL({\cal P}_n)-\alpha n^{(d-p)/d} \le \left\{ \ba{ll}
Cn^{(d-1-p)/d}&\mb{ for $0<p<d-1$}\\
C(\log n \vee 1)&\mb{ for $p=d-1\neq 1$}\\
C&\mb{ for $d-1<p<d$} \\
C&\mb{ for $p=d-1=1$}. \ea \right.
$$
where $\al:=\al(L,d,p)$ is a finite constant given by \eq{1.6}.
\ignore{ and
$$
\alpha n^{(d-p)/d}-C \leq EL({\cal P}_n,[0,1]^d,p).
$$
If also there exists a superadditive Euclidean functional $L^*$ of
power $p$ such that $(L,L^*)$ is an Euclidean functional pair of
power $p$ and $L$ is a strong Euclidean functional, then \be
EL({\cal P}_n,[0,1]^d,p)\le \alpha n^{(d-p)/d}+\left\{ \ba{ll}
Cn^{(d-1-p)/d}&\mb{ for $0<p<d-1$}\\
C(\log n \vee 1)&\mb{ for $p=d-1$}\\
C&\mb{ for $d-1<p<d$}. \ea \right.\nn \ee }

\noin {\bf Proof.} By the usual subadditive argument (see page 54
of \cite{Y98}) for $L$, \ignore{ In the proof of Theorem 1.1 (See
\cite{S88}), we see
$$
\lim_{t \ra \infty}\frac{\varphi (t)}{t^d}=\lim_{t \ra
\infty}\frac{\varphi^* (t)}{t^d}=\alpha.
$$
As we have mentioned,
$$
\varphi (t)=EL({\cal P},[0,t]^d,p)=t^pEL({\cal
P}_{t^d},[0,1]^d,p).
$$
Thus
$$
\lim_{t \ra \infty}\frac{EL({\cal P}_{t^d})}{t^{d-p}}=\lim_{n \ra
\infty}\frac{EL({\cal P}_{n})}{n^{(d-p)/d}}=\alpha.
$$
By the subadditivity condition \eq{1.6},
$$
L({\cal P},[0,t]^d,p) \leq \sum_{j=1}^{m^d} L({\cal P} \cap
Q_i,Q_i,p) +Cm^{d-p}t^p.
$$
Take the expectation with $t=mu$
$$
\frac{\varphi (mu)}{(mu)^d} \leq \frac{\varphi (u)}{u^d}
+Cu^{-(d-p)}.
$$
Letting $m \ra \infty$ with $n=u^d$ yields
$$
\alpha \leq \frac{\varphi (u)}{u^d} +Cu^{-(d-p)}.
$$
$$
\alpha  \leq \frac{EL({\cal
P}_{u^d}[0,1]^d,p)}{u^{d-p}}+Cu^{-(d-p)}=\frac{EL({\cal
P}_{n},[0,1]^d,p)}{n^{(d-p)/d}} +Cn^{-(d-p)/d}.
$$
Thus, }
$$
-C \leq EL({\cal P}_{n})-\alpha n^{(d-p)/d}.
$$
By the same argument for $L^*$,
$$
EL^*({\cal P}_{n})-\alpha n^{(d-p)/d} \le C.
$$
Now, use \eq{1.11} to get the Lemma. \ignore{ Suppose now that
$(L,L^*)$ is  an Euclidean functional pair of power $p$, then
$$
L^*({\cal P},[0,t]^d,p) \geq \sum_{j=1}^{m^d} L^*({\cal P} \cap
Q_i,Q_i,p) - C m^{d-p}t^{p}.
$$
Take the expectation with $t=mu$
$$
\frac{\varphi^* (mu)}{(mu)^d} \geq \frac{\varphi^* (u)}{u^d} - C
u^{-(d-p)}.
$$
Letting $m \ra \infty$ with $n=u^d$ yields
$$
\alpha \geq \frac{\varphi^* (u)}{u^d}- C u^{-(d-p)}.
$$
$$
\varphi^* (u)=EL^*({\cal P},[0,u]^d,p)=u^pEL^*({\cal
P}_{u^d},[0,1]^d,p).
$$
$$
\alpha \geq \frac{EL^*({\cal P}_{u^d},[0,1]^d,p)}{u^{d-p}}- C
u^{-(d-p)} =\frac{EL^*({\cal P}_{n},[0,1]^d,p)}{n^{(d-p)/d}}- C
n^{-(d-p)/d}.
$$
$$
\alpha n^{(d-p)/d} + C \geq EL^*({\cal P}_{n},[0,1]^d,p).
$$
Since $L$ and $L^*$ are close in mean, \bea EL({\cal
P}_{n},[0,1]^d,p) &\le& EL^*({\cal P}_{n},[0,1]^d,p)+\left\{
\ba{ll}
Cn^{(d-1-p)/d}&\mb{ for $0<p<d-1$}\\
C\log n&\mb{ for $p=d-1$}\\
C&\mb{ for $d-1<p<d$} \ea
\right.\nn\\
&\le& \alpha n^{(d-p)/d}+\left\{ \ba{ll}
Cn^{(d-1-p)/d}&\mb{ for $0<p<d-1$}\\
C(\log n \vee 1)&\mb{ for $p=d-1$}\\
C&\mb{ for $d-1<p<d$}. \ea \right. \nn \eea This concludes the
Lemma 2.1.}
\qed\\

\noin {\bf Lemma 2.} Let $N_n$ is a Poisson random variable with
mean $n$ which is independent of i.i.d.\ uniform points
$\{X_i; 1 \le i <\infty\}$ on $[0,1]^d$.

(i) If  $L$ is a strong Euclidean functional of power $0<p<d$,
then \bea \left|EL({\cal U}_n)-EL({\cal U}_{N_n})\right| \leq
C(n^{(d-1-p)/2(d-1)} \vee 1). \lbl{2.1} \eea

(ii) If  $L$ is a very strong Euclidean functional of power
$0<p<d$, then \be \left|EL({\cal U}_n)-EL({\cal U}_{N_n})\right|
\le Cn^{1/2 - p/d}. \lbl{2.2} \ee

\noin {\bf Proof.} (i) If  $L$ is a strong Euclidean functional of
power $0<p<d-1$, by \eq{1.12}, \eq{a} and by Jensen's inequality
we have \bea \left|EL({\cal U}_n)-EL({\cal U}_{N_n})\right| &\le&
\left|E(L({\cal U}_n) - L({\cal U}_{N_n})) {\bold 1}_{\{|N_n-n|\leq n^{(d-1)/d}\}}\right|  \nn \\
&&
+\left|E(L({\cal U}_n) - L({\cal U}_{N_n})) {\bold 1}_{\{n^{(d-1)/d}\leq |N_n-n|\leq n/2\}}\right|  \nn \\
&   & + \left|E(L({\cal U}_n) - L({\cal U}_{N_n})) {\bold
1}_{\{N_n <
n/2\}}\right|\nn\\
&&
+ \left|E(L({\cal U}_n) - L({\cal U}_{N_n})) {\bold 1}_{\{N_n > 3n/2\}}\right| \nn \\
&\leq& C E(|N_n-n|^{(d-1-p)/(d-1)}\vee 1)
+ C n^{-p/d} E|N_n - n|   \nn \\
&   & + C n^{(d-p)/d} P(N_n < n/2) + C EN_n^{(d-p)/d}
{\bold 1}_{\{N_n > 3n/2\}} \nn\\
&\leq& C(n^{(d-1-p)/2(d-1)} \vee 1)
+ Cn^{1/2 - p/d}   \nn \\
&   & + C n^{(d-p)/d} P(N_n < n/2) + C EN_n^{(d-p)/d} {\bold
1}_{\{N_n > 3n/2\}}. \nn \eea Since $N_n$ has a very light tail,
the last two terms of the left hand side in the above inequality
are negligible and \eq{2.1} follows. The argument for $d-1 \leq p
< d$ is same.

(ii) If  $L$ is a very strong Euclidean functional of power
$0<p<d$, there is a standard argument from Theorem 5.2 of Yukich
(1998): By \eq{a} we have \bea \left|EL({\cal U}_n)-EL({\cal
U}_{N_n})\right| &\le&
C n^{-p/d} E\left|N_n - n \right|{\bold 1}_{\{1\leq |N_n-n|\leq n/2\}}  \nn\\
&&+ C n^{(d-p)/d} P(N_n < n/2) + C EN_n^{(d-p)/d} {\bold 1}_{\{N_n
> 3n/2\}}. \nn \eea Thus, \eq{2.2} follows.
\qed\\

\noin{\bf Proof of Theorem 1.}
Theorem 1 follows from Lemmas 1 and 2. \qed\\

\noin {\bf Lemma 3.} For the superadditive Euclidean functional
$L^*$ of MM, MST, TSP with $0<p<d$ and  SMST, RSMST with $0<p\le
1$, let $N_B({\cal P}_n)$ be the number of points in  ${\cal P}_n$
which are connected to  the boundary $\partial [0,1]^d$  in the
minimal graph which is used to calculate $L^*$. Using the same
minimal graph let $L_B({\cal P}_n)$ be the sum of the $p$-th
power-weighted lengths of the edges connecting points in ${\cal
P}_n$ to the boundary $\partial [0,1]^d$. Then,
$$
EN_B({\cal P}_n) \le Cn^{(d-1)/d},
$$
$$
EL_B({\cal P}_n) \le \left\{ \ba{ll}
Cn^{(d-1-p)/d}&\mb{ for $0<p<d-1$}\\
C(\log n \vee 1)&\mb{ for $p=d-1\neq 1$}\\
C&\mb{ for $d-1<p<d$}\\
C&\mb{ for $p=d-1=1$}. \ea \right.
$$
The same estimates for ${\cal U}_n$ instead of ${\cal P}_n$ also
holds.

 \noin {\bf Proof.} We just follow the
argument for Lemma 3.8 of Yukich (1998) or Lemma 4 of Lee (1999)
and dig out the quantities of interest. We skip its proof.
\ignore{
 Since the case $1 \le p <d$ has been known quite sometime,
we prove the case  $0<p<1$ only. Furthermore, here we prove the
Lemma for the MM and leave all the other cases to the reader as an
exercise.

The idea is that, if two points are close to each other, then at
most one of these two points is matched to the boundary points in
the optimal dual matching since matching these two points to each
other is more economic then matching these two points to the
boundary points. Let $D_j$, $j=0,1,\ldots,k$, be a box of the form
$ D_j =
[1/2-1/6-\sum_{l=1}^j3^{-1}2^{-l},1/2+1/6+\sum_{l=1}^j3^{-1}2^{-l}]^d.
$ We choose $k$ so that the moat $[0,1]^d \setminus \cup_{j=0}^k
D_j$ has a width of order $n^{-1/d}$, i.e., $ 3^{-1}2^{-k} \geq
n^{-1/d} > 3^{-1}2^{-k-1}. $ In order to quantify the closeness of
the two points we partition $\cup_{j=0}^k D_j$.  First we
partition $D_{j+1} \setminus D_j$, with a convention
$D_{-1}=\emptyset$, into subboxes of edge length $3^{-1}2^{-j-1}$.
We part ion these subboxes further into subboxes of edge length
$3^{-1}2^{-j-m-1}$ where $m$ is the smallest positive integer such
that $1 > (2^{-m}\sqrt{d})^p$. Suppose that there are two points
in the same subbox.  Matching these two points to each other costs
at most $(3^{-1}2^{-j-m-1}\sqrt{d})^p$ whereas matching these two
points to the boundary points costs at least $(3^{-1}2^{-j-1})^p$.
By the choice of $m$ we have $(3^{-1}2^{-j-m-1}\sqrt{d})^p <
(3^{-1}2^{-j-1})^p$, i.e., since matching these two points each
other is more economic then matching these two points to the
boundary points, at most one of these two points is matched to the
boundary points in the optimal dual matching.

Therefore, we can decompose $N_B({\cal P}_n)$ into two parts;  the
number $N_m({\cal P}_n)$ of vertices in  the moat $[0,1]^d
\setminus \cup_{j=0}^kD_j$
 which is connected to  the boundary of
$[0,1]^d$  in the optimal dual graph and
  the number $N_p({\cal P}_n)$  of vertices in   the partition
$\cup_{j=0}^kD_j$
 which is connected to  the boundary of
$[0,1]^d$  in the optimal dual graph.
$$
EN_m({\cal P}_n) \le Cn^{(d-1)/d}.
$$
$$
N_p({\cal P}_n) \le C\sum_{j=0}^{k} 2^{j(d-1)}\le Cn^{(d-1)/d}.
$$
Now, the Lemma for $N_B({\cal P}_n)$ follows.

Also, we can decompose $L_B({\cal P}_n)$ into two parts;  the sum
$L_m({\cal P}_n)$ of the $p$-th power of the lengths of the edges
connecting vertices in the moat $[0,1]^d \setminus
\cup_{j=0}^kD_j$ with the boundary of $[0,1]^d$ and
  the sum $L_p({\cal P}_n)$ of the $p$-th power
of the lengths of the edges connecting vertices in the partition
$\cup_{j=0}^kD_j$ with the boundary of $[0,1]^d$.
$$
EL_m({\cal P}_n) \le Cn^{-p/d}EN_m({\cal P}_n) \le Cn^{(d-1-p)/d}.
$$
$$
L_p({\cal P}_n) \le C\sum_{j=0}^{k} 2^{-jp} 2^{j(d-1)}\le
Cn^{(d-1-p)/d}.
$$
Now, the Lemma for $L_B({\cal P}_n)$ follows }
\qed\\

\noin {\bf Lemma 4.} For the subadditive Euclidean functional $L$
of MM, MST, TSP with $0<p<d$ and  SMST, RSMST with $0<p\le 1$, \be
0 \le EL({\cal P}_n)-EL^*({\cal P}_n)\leq \left\{ \ba{ll}
Cn^{(d-1-p)/d}&\mb{ for $0<p<d-1$}\\
C(\log n \vee 1)&\mb{ for $p=d-1\neq 1$}\\
C&\mb{ for $d-1<p<d$} \\
C&\mb{ for $p=d-1=1$}. \ea \right. \lbl{2.3} \ee The same
estimates for ${\cal U}_n$ instead of ${\cal P}_n$ also holds.\\
\noin {\bf Proof.} Since the case $1 \le p <d$ has been known
quite sometime (see Lemma 3.10 of Yukich (1998)), we consider the
case $0<p<1$ only. Even in this case, the argument for the case $1
\le p <d-1$ of Lemma 3.10 of Yukich (1998) works; we just need to
use the corresponding estimates for $0<p<1$ in Lemma 3.
\qed\\

\noin {\bf Lemma 5.} For the subadditive Euclidean functional $L$
of MM, MST, TSP with $0<p<d$ and  SMST, RSMST with $0<p\le 1$, $L$
satisfies \eq{1.12}. \ignore{ the following inequality. \bea
& & \Big|EL({\cal U}_{n \pm k},[0,1]^d,p)-EL({\cal U}_n,[0,1]^d,p)\Big| \nn \\
&=& \left\{ \ba{lll}
C(k^{(d-1-p)/(d-1)}\vee 1)&\mb{ for $1 \le k \le n^{(d-1)/d}$}\\
Ckn^{-p/d} &\mb{ for $n^{(d-1)/d} < k \le n/2$, $0<p<d-1$}\\
C &\mb{ for $n^{(d-1)/d} < k \le n/2$, $d-1\leq p<d$} \ea \right.
\nn \eea }

\noin {\bf Proof.} Since all the arguments are similar, here we
prove the Lemma for the MST case only. We leave all the other
cases to the reader as an exercise. First we claim that  for $1
\le k \le n^{(d-1)/d}$ \be EL({\cal U}_{n + k}) \le EL({\cal
U}_n)+ C(k^{(d-1-p)/(d-1)}\vee 1). \lbl{2.4} \ee Let
$\{X_1,\ldots,X_{n+k} \}$ be the $n+k$ i.i.d.\ uniform points on
$[0,1]^d$. By renaming these points $X_i:=(X_i^1,\ldots,X_i^d)$ we
may assume that the first coordinates of these $n+k$  points are
increasing, i.e., $ X_{1}^{1} < X_{2}^1 < \cdots < X_{n+k}^1 $.
With this renaming, let $ \tau := 1-X_{n+1}^1 $. After
constructing  an MST on $\{ X_1,\ldots,X_n\}$ and another MST on
$\{X_{n+1},\ldots,X_{n+k}\}$, by adding an edge $(X_n,X_{n+1})$ we
have a spanning tree on ${\cal U}_{n+k}$. So, we have  $$ EL({\cal
U}_{n+k}) \leq EL(\{ X_1,\ldots,X_n\}) +
EL(\{X_{n+1},\ldots,X_{n+k}\}) + C.$$ Let $X_j^* := ((1-\tau)^{-1}
X_j^1,X_j^2, X_j^3, \ldots, X_j^d), 1 \leq j \leq n$. Then, since
$\{ X_1^*,\ldots, X_n^* \}$ are $n$ i.i.d.\ uniform samples
from $[0,1]^d$,
$$
EL(\{X_1,\ldots,X_n\}) \leq EL(\{X_1^*,\ldots, X_n^* \}) =
EL({\cal U}_n),
$$
and hence \be EL({\cal U}_{n+k}) \leq EL({\cal U}_n) +
EL(\{X_{n+1},\ldots,X_{n+k}\}) + C. \lbl{2.5} \ee

Now, we estimate $EL(\{X_{n+1},\ldots,X_{n+k}\})$. Let $Y_i$, $n+1
\leq i \leq n+k$, be the projection of $X_i$ to $\{1\} \times
[0,1]^{d-1}$. We construct an MST on $\{Y_{n+1},\ldots,Y_{n+k}\}$
and then using this MST we construct a spanning tree on
$\{X_{n+1},\ldots,X_{n+k}\}$ by adding an edge $(X_i,X_j)$ if an
edge  $(Y_i,Y_j)$ is in this MST (see Figure 1). For this edge
$(X_i,X_j)$
$$
|X_i-X_j|^p \le \left(|Y_i-Y_j|+2\tau\right)^p \le
C|Y_i-Y_j|^p+C\tau^p.
$$

\setlength{\unitlength}{1.0mm} 
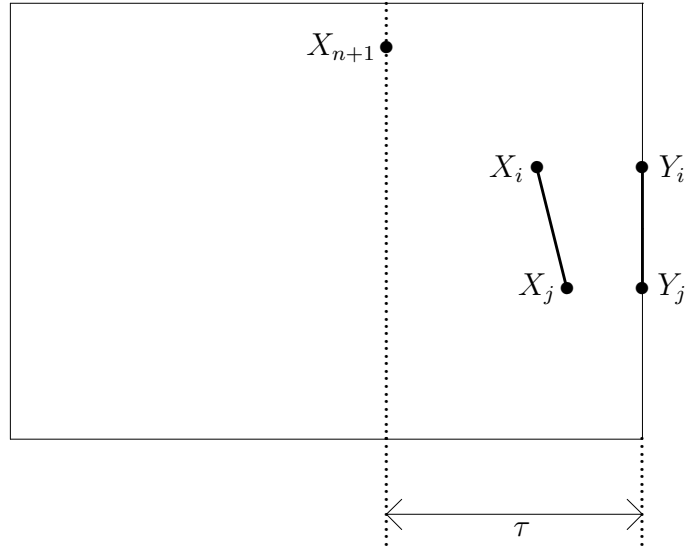
\begin{figure}
\centering
\begin{picture}(90,80)
\thinlines \drawpath{4.0}{76.0}{88.0}{76.0}
\drawpath{88.0}{76.0}{88.0}{76.0} \drawpath{4.0}{76.0}{4.0}{18.0}
\drawpath{4.0}{18.0}{88.0}{18.0}
\drawdotline{54.0}{76.0}{54.0}{18.0}
\drawpath{88.0}{76.0}{88.0}{46.0}
\drawpath{88.0}{46.0}{88.0}{46.0}
\drawpath{88.0}{38.0}{88.0}{18.0}
\drawcenteredtext{54.0}{70.0}{$\bullet$}
\drawcenteredtext{88.0}{54.0}{$\bullet$}
\drawcenteredtext{88.0}{38.0}{$\bullet$}
\drawcenteredtext{74.0}{54.0}{$\bullet$}
\drawcenteredtext{78.0}{38.0}{$\bullet$} \Thicklines
\drawpath{74.0}{54.0}{78.0}{38.0}
\drawpath{88.0}{54.0}{88.0}{38.0} \thinlines
\drawdotline{54.0}{18.0}{54.0}{4.0}
\drawdotline{88.0}{18.0}{88.0}{4.0}
\drawpath{54.0}{8.0}{88.0}{8.0} \drawpath{54.0}{8.0}{56.0}{10.0}
\drawpath{54.0}{8.0}{56.0}{6.0} \drawpath{86.0}{10.0}{88.0}{8.0}
\drawpath{88.0}{8.0}{86.0}{6.0}
\drawcenteredtext{72.0}{6.0}{$\tau$}
\drawcenteredtext{36.0}{52.0}{}
\drawcenteredtext{70.0}{54.0}{$X_i$}
\drawcenteredtext{74.0}{38.0}{$X_j$}
\drawcenteredtext{92.0}{54.0}{$Y_i$}
\drawcenteredtext{92.0}{38.0}{$Y_j$}
\drawcenteredtext{48.0}{70.0}{$X_{n+1}$}
\end{picture}

\caption[short caption here]{A way to construct a spanning tree on
$\{X_{n+1},\ldots,X_{n+k} \}$; add an edge $(X_i,X_j)$ if an edge
$(Y_i,Y_j)$ is in the MST on $\{Y_{n+1},\ldots,Y_{n+k} \}$.}
\label{fig:test} 
\end{figure}

Since $E\tau^p \leq C(k/n)^p$, we have then by \eq{*} \bea
EL(\{X_{n+1},\ldots,X_{n+k} \}) &\leq&
CEL(\{Y_{n+1},\ldots,Y_{n+k} \},[0,1]^{d-1} \times \{1\},p) + C k
E\tau^p\nn\\
&\le& C (k^{(d-1-p)/(d-1)}\vee
1) + C k E\tau^p \nn\\
&\leq& C (k^{(d-1-p)/(d-1)} \vee 1).\lbl{2.6} \eea Therefore,
\eq{2.4} follows from \eq{2.5}-\eq{2.6}.

Second,  for $n^{(d-1)/d} < k \leq n/2$ we iterate the above
argument $[k/n^{(d-1)/d}]+1$ times and we have \bea EL({\cal
U}_{n+ k}) &\leq&
EL({\cal U}_n) + C([k/n^{(d-1)/d}] +1)(n^{(d-1-p)/d} \vee 1)\nn \\
&\leq& EL({\cal U}_n) + \left\{ \ba{ll}
Ckn^{-p/d} &\mb{ for $0<p<d-1$}\\
Ckn^{-(d-1)/d} &\mb{ for $d-1\leq p<d$}. \ea \right. \lbl{2.7}
\eea

\setlength{\unitlength}{1.0mm} 
\begin{figure}
\centering
\begin{picture}(90,80)
\thinlines \drawpath{4.0}{76.0}{80.0}{76.0}
\drawpath{80.0}{76.0}{80.0}{76.0} \drawpath{4.0}{76.0}{4.0}{18.0}
\drawpath{4.0}{18.0}{4.0}{18.0} \drawpath{80.0}{76.0}{80.0}{18.0}
\drawpath{80.0}{18.0}{4.0}{18.0} \drawpath{28.0}{76.0}{26.0}{76.0}
\drawdotline{30.0}{76.0}{30.0}{76.0}
\drawdotline{34.0}{76.0}{34.0}{18.0}
\drawdotline{64.0}{76.0}{64.0}{76.0}
\drawdotline{60.0}{76.0}{60.0}{18.0}
\drawpath{22.0}{40.0}{70.0}{40.0}
\drawpath{22.0}{40.0}{22.0}{40.0}
\drawdotline{22.0}{40.0}{38.0}{30.0}
\drawpath{38.0}{30.0}{56.0}{26.0}
\drawdotline{70.0}{40.0}{56.0}{26.0}
\drawdotline{34.0}{18.0}{34.0}{4.0}
\drawdotline{60.0}{18.0}{60.0}{4.0}
\drawdotline{60.0}{4.0}{60.0}{4.0}
\drawpath{34.0}{10.0}{60.0}{10.0}
\drawpath{34.0}{10.0}{36.0}{12.0}
\drawpath{36.0}{12.0}{36.0}{12.0} \drawpath{34.0}{10.0}{36.0}{8.0}
\drawpath{60.0}{10.0}{58.0}{12.0}
\drawpath{58.0}{12.0}{58.0}{12.0} \drawpath{60.0}{10.0}{58.0}{8.0}
\drawpath{58.0}{8.0}{58.0}{8.0}
\drawcenteredtext{34.0}{50.0}{$\bullet$}
\drawcenteredtext{60.0}{60.0}{$\bullet$}
\drawcenteredtext{22.0}{40.0}{$\bullet$}
\drawcenteredtext{70.0}{40.0}{$\bullet$}
\drawcenteredtext{38.0}{30.0}{$\bullet$}
\drawcenteredtext{56.0}{26.0}{$\bullet$}
\drawcenteredtext{44.0}{8.0}{$\tau_i$}
\drawcenteredtext{30.0}{50.0}{$X_i$}
\drawcenteredtext{66.0}{60.0}{$X_{i+k}$}
\drawcenteredtext{18.0}{40.0}{$X_l$}
\drawcenteredtext{74.0}{40.0}{$X_{l'}$}
\drawcenteredtext{38.0}{28.0}{$X_j$}
\drawcenteredtext{56.0}{24.0}{$X_{j'}$}
\end{picture}

\caption[short caption here]{A way to construct a connected graph
on ${\cal B}_i$; add an edge $(X_l,X_{l^{'}})$, $X_l \in {\cal
B}_i$, $X_{l^{'}} \in {\cal B}_i$, if $(X_l,X_j)$, $X_j \in {\cal
W}_i$, and $(X_{l^{'}},X_{j^{'}})$, $X_{j^{'}} \in {\cal W}_i$,
were edges in the original MST $T({\cal U}_n)$ and if
$(X_j,X_{j^{'}})$ is an edge in the MST $T({\cal W}_i)$.}
\label{fig:test} 
\end{figure}
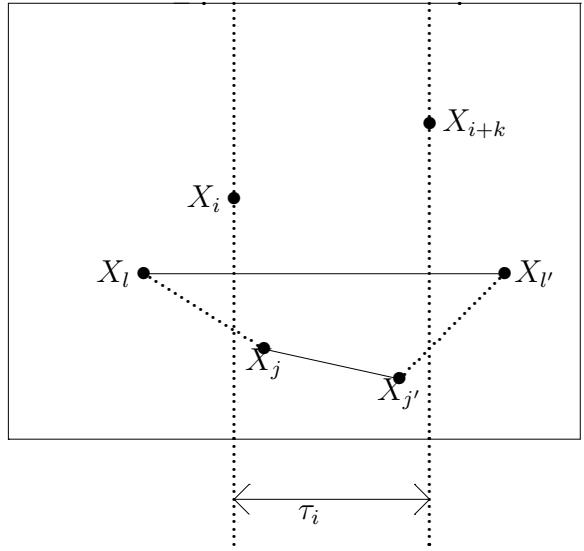

Third, we claim that  for $1 \le k \le n^{(d-1)/d}$ \be EL({\cal
U}_{n - k}) \le EL({\cal U}_n)+ C(k^{(d-1-p)/(d-1)} \vee
1).\lbl{2.8} \ee As we did before, by renaming the $n$ random
points ${\cal U}_n:= \{ X_1,\ldots,X_n\}$ we may assume that the
first coordinates of these $n$  points are increasing, i.e., $
X_{1}^{1} < X_{2}^1 < \cdots < X_{n}^1 $. For $1\le i\leq n-k+1$,
let ${\cal B}_i = {\cal U}_n \setminus \{ X_i,\ldots,X_{i+k-1} \}$
and let $\tau_i := X_{i+k}^1 - X_i^1$. We construct an MST
$T({\cal U}_n)$ on ${\cal U}_n$. From the MST $T({\cal U}_n)$, we
remove all the points $X_j$, $i \leq j \leq i+k-1$, and all the
edges in the MST $T({\cal U}_n)$ which use those $X_j$ as one end.
The resulting graph $D_i$ on ${\cal B}_i$ is disconnected. By
adding some extra edges to $D_i$ we will construct a connected
graph on ${\cal B}_i$. To do this we first collect those $X_j$, $i
\leq j \leq i+k-1$, for which $X_j$ was connected to a remainder
point $X_k \in {\cal B}_i$ by an edge in the original MST $T({\cal
U}_n)$. We call the set of all those collected $X_j$ as ${\cal
W}_i$. We construct an MST $T({\cal W}_i)$ on ${\cal W}_i$. Using
this MST $T({\cal W}_i)$ we have the following addition rule of
edges to the disconnected graph $D_i$ on ${\cal B}_i$; we add an
edge $(X_l,X_{l^{'}})$, $X_l \in {\cal B}_i$, $X_{l^{'}} \in {\cal
B}_i$, if $(X_l,X_j)$, $X_j \in {\cal W}_i$, and
$(X_{l^{'}},X_{j^{'}})$, $X_{j^{'}} \in {\cal W}_i$, were edges in
the original MST $T({\cal U}_n)$ and if $(X_j,X_{j^{'}})$ is an
edge in the MST $T({\cal W}_i)$. By adding these extra edges to
$D_i$ we get a connected graph on ${\cal B}_i$ (see Figure 2).

In this case
$$
|X_l-X_{l'}|^p \le
C|X_l-X_{j}|^p+C|X_j-X_{j'}|^p+C|X_{j'}-X_{l'}|^p.
$$
Since the degree of a vertex in an MST on the set of points in
$\mathR^d$ is bounded by a universal constant which depends only
on the dimension $d$, by the argument of \eq{2.6} \bea &&
EL({\cal B}_i)\nn\\
&\leq& EL({\cal U}_n) + C\sum_{j=i}^{i+k-1} \sum_{(X_j,X_k) \in
T({\cal U}_n)}E|X_j-X_k|^p + C EL_{}({\cal
W}_i)\nn\\
&\leq& EL({\cal U}_n) + C\sum_{j=i}^{i+k-1} \sum_{(X_j,X_k) \in
T({\cal U}_n)}E|X_j-X_k|^p + C(k^{(d-1-p)/(d-1)}\vee 1).\lbl{2.9}
\eea Let $X_j^*:= X_j - \tau_i e_1$ for $i+k \leq j \leq n,$ where
$e_1:=(1,0,\ldots,0))$ is the first unit vector, and let ${\cal
B}_i^* := \{X_1,\ldots,X_{i-1},X_{i+k}^*,\ldots,X_n^*\}$.  ${\cal
B}_i^*$ then consists of the $n-k$ i.i.d.\ uniform samples from
$[0,1]^{d-1} \times [0,1-\tau_i]$. By its construction,
$$
L_{}({\cal B}_i^*) \leq L_{}({\cal B}_i).
$$
Multiplying the first coordinates of the points of ${\cal B}_i^*$
by $(1-\tau_i)^{-1}$ we have $n$ i.i.d.\ uniform samples from
$[0,1]^d$. Hence, since $(1-\tau_i)^{-p} \le 1+C\tau_i$ for the
case $\tau_i \leq 1/2$, and since $1+C\tau_i$ is increasing in
$\tau_i$ and $E(L_{}({\cal B}_i^*)|\tau_i)$ is decreasing in
$\tau_i$, \bea &&
EL_{}({\cal U}_{n-k})\nn\\
&\leq& E\left((1-\tau_i)^{-p} E(L_{}({\cal B}_i^*)|\tau_i)|\tau_i
\leq 1/2\right)P(\tau_i \leq 1/2) + E\left(L_{}({\cal
U}_{n-k});\tau_i>1/2\right)
\nn\\
&\leq& E\left((1+C\tau_i) E(L_{}({\cal B}_i^*)|\tau_i)|\tau_i \leq
1/2\right)P(\tau_i \leq 1/2)+
e^{-Cn}\nn\\
&\leq& \left(E\left(1+C\tau_i|\tau_i \leq 1/2\right)
E\left(E(L_{}({\cal B}_i^*)|\tau_i)|\tau_i \leq
1/2\right)\right)P(\tau_i \leq 1/2)+
e^{-Cn}\nn\\
&\leq& \left(1+C\frac{k}{n+1}\right) E\left(E(L_{}({\cal
B}_i^*)|\tau_i)|\tau_i \leq 1/2\right)+
e^{-Cn}\nn\\
&\leq& \left(1+C\frac{k}{n+1}\right) \frac{EL_{}({\cal
B}_i^*)}{P(\tau_i \leq 1/2)}+
e^{-Cn}\nn\\
&\leq& \left(1+C\frac{k}{n+1}\right) \frac{EL_{}({\cal
B}_i)}{P(\tau_i \leq 1/2)}+
e^{-Cn} \mb{ (by \eq{2.9})}\nn\\
&\le& \left(1+C\frac{k}{n}\right)EL({\cal U}_n) +
C\sum_{j=i}^{i+k-1} \sum_{(X_j,X_k) \in T({\cal U}_n)}E|X_j-X_k|^p
+ C(k^{(d-1-p)/(d-1)}\vee 1). \nn \eea Since $k \leq n^{(d-1)/d}$,
we have
$$
\frac{k}{n}EL({\cal U}_n) \leq C \frac{k}{n} n^{(d-p)/d} \leq
Ck^{(d-1-p)/(d-1)}.
$$
So, by combining the terms $C \frac{k}{n} EL({\cal U}_n)$ and $C
(k^{(d-1-p)/(d-1)} \vee 1)$ and by increasing the constant $C$ in
the term $C(k^{(d-1-p)/(d-1)} \vee 1)$ we have \be EL({\cal
U}_{n-k}) \le EL({\cal U}_n) + C\sum_{j=i}^{i+k-1} \sum_{(X_j,X_k)
\in T({\cal U}_n)}E|X_j-X_k|^p + C(k^{(d-1-p)/(d-1)}\vee 1).
\lbl{2.10} \ee \ignore{ If $\tau_i \leq 1/2$, then \bea
(1-\tau_i)^{-p} \leq \left\{ \ba{ll}
1+2\tau_i &\mb{ for $0<p\le 1$}\\
1+2^{p+1}\tau_i &\mb{ for $1< p<d$}. \ea \right. \lbl{2.23} \eea
Since two cases are proved similarly, we will show that the only
for $1<p<d$. Hence averaging \eq{2.19} over the event $[\tau_i
\leq 1/2]$ gives
$$
EL_{}({\cal U}_{n-k}) \leq E[(1+2^{p+1}\tau_i) E(L_{}({\cal
B}_i^*)|\tau_i)]/P[\tau_i \leq 1/2]
$$
} Here comes the highlight of the argument; we use the symmetry
argument to estimate the term $\sum_{j=i}^{i+k-1} \sum_{(X_j,X_k)
\in T({\cal U}_n)}E\left|X_j-X_k\right|^p$. Since
$$
\sum_{i=1}^{n-k+1}\sum_{j=i}^{i+k-1} \sum_{(X_j,X_k) \in T({\cal
U}_n)}E\left|X_j-X_k\right|^p \le CkEL_{}({\cal U}_n),
$$
averaging \eq{2.10} over $1\le i \le n-k+1$ (by \eq{a}) we have
\eq{2.8}.

Last,  for $n^{(d-1)/d} < k \leq n/2$ we iterate the above
argument $\left[k/n^{(d-1)/d}\right]+1$ times and we have \bea
EL({\cal U}_{n-k}) &\leq&
EL({\cal U}_n) + C([k/n^{(d-1)/d}] +1) (n^{(d-1-p)/d} \vee 1)\nn \\
&\leq& EL({\cal U}_n) + \left\{ \ba{ll}
Ckn^{-p/d} &\mb{ for $0<p<d-1$}\\
Ckn^{-(d-1)/d} &\mb{ for $d-1\leq p<d$}. \ea \right. \lbl{2.11}
\eea
\qed\\

For the Euclidean functional $L$ of MST with $0<p<d$, $L$
satisfies the add-one bound \eq{1.13} as shown by Redmond and
Yukich (1994). In Lemma 6 below we show that for the case $0<p\leq
1$ many typical $L$ also satisfy the add-one bound \eq{1.13} and
hence we provide some affirmative answers to the issue raised in
p.\ 55 of Yukich (1998). However, we cannot handle the case
$1<p<d$ for those $L$ and more seriously we cannot prove the
add-one bound for the minimal matching Euclidean functional
$L_{MM}$. So, we think the add-one bound \eq{1.13} condition is
very restrictive.\\

\noin{\bf Lemma 6.}
If $L$ is either the MST Euclidean functional $L$ with
$0<p<d$ or the TSP, SMST, RSMST  Euclidean functional  with $0<p\le 1$,  then $L$ satisfies
\eq{1.13}.

\noin {\bf Proof.} Since all the arguments are similar, here we
prove the Lemma for the TSP case only. We leave all the other
cases to the reader as an exercise.

Fix $0<p\le 1$. First we claim that \be
EL({\cal U}_{n}) \leq EL({\cal U}_{n+1})+ Cn^{-p/d}. \lbl{2.12}
\ee There are two edges $(X_j,X_i)$ and $(X_i,X_k)$ adjacent to
$X_i$ in the minimal tour $T({\cal U}_{n+1})$. Remove these two
and add the edge $(X_j,X_k)$. Then, we have a tour on ${\cal
U}_{n+1}\setminus \{X_i\}$. Since
$$
\left|X_j-X_k\right|^p \le
\left(\left|X_j-X_i\right|+\left|X_i-X_k\right|\right)^p \le
C\left|X_j-X_i\right|^p+C\left|X_i-X_k\right|^p,
$$
we have
$$
EL({\cal U}_{n+1}\setminus \{X_i\}) \leq EL({\cal U}_{n+1}) + C
\sum_{(X_i,X_j) \in T({\cal U}_{n+1})}E \left| X_i - X_j\right|^p.
$$
Since $EL({\cal U}_{n+1}\setminus \{X_i\})=EL({\cal U}_n)$, and
since  $\sum_{i=1}^{n+1}\sum_{(X_i,X_j) \in T({\cal U}_{n+1})}E
\left| X_i - X_j\right|^p = 2EL({\cal U}_{n+1})\le Cn^{(d-p)/d}$
(by \eq{a}), averaging the above inequality over $1\leq i \leq
n+1$ we have \eq{2.12}. \ignore{
$$
\sum_{i=1}^{n} L({\cal U}_{n},\hat{X_i},[0,1]^d) \leq
\sum_{i=1}^{n} L({\cal U}_{n+1},[0,1]^d,p) + C_{14} \sum_{i=1}^{n}
\sum_{j \in N(i)} \Big| X_i - X_j \Big|^p.
$$
Taking expectations and noting that the expectation of the double
sum is bounded by $Cn^{(d-p)/d}$ gives
$$
nEL({\cal U}_{n},[0,1]^d,p) \leq n EL({\cal U}_{n+1},[0,1]^d,p) +
Cn^{(d-p)/d}.
$$
Dividing by $n$
$$
EL({\cal U}_{n},[0,1]^d,p) \leq EL({\cal U}_{n+1},[0,1]^d,p) +
Cn^{-p/d}.
$$
}

Now, we claim that \be EL({\cal U}_{n+1}) \leq EL({\cal U}_{n})+
Cn^{-p/d}.\lbl{2.13} \ee Let $T({\cal U}_n)$ be the minimal tour
on ${\cal U}_n$. For a given $X_{n+1}$, among $n$ i.i.d.\
uniform points ${\cal U}_n$ find a nearest point $X_i$ to
$X_{n+1}$ and let $(X_i,X_j)$ be an edge in the minimal tour
$T({\cal U}_n)$. By removing this edge and by adding
$(X_i,X_{n+1})$ and $(X_{n+1},X_j)$ we construct a tour on ${\cal
U}_{n+1}$.
Since $0<p\le 1$, we have $A^p\le B^p+C^p$ for a triangle with side lengths $A$, $B$, $C$.
So, using this triangular inequality
 \bea L({\cal U}_{n+1}) &\leq&
\L({\cal U}_{n}) - \left|X_i - X_{j}\right|^p + \left|X_i - X_{n+1}\right|^p + \left|X_{j} - X_{n+1}\right|^p \nn \\
&\leq&
L({\cal U}_n) + 2\left|X_i - X_{n+1}\right|^p \nn \\
&=&   L({\cal U}_n) + 2 \min_{1 \le k \leq n} \left|X_k
-X_{n+1}\right|^p. \nn \eea By taking  expectations in the above
inequality we have \eq{2.13}.
\qed\\

\noin{\bf Proof of Theorem 2.} (i) follows from Lemmas 4 and 5
and (ii) follows from Lemmas 4, 5, and 6. \qed\\

\section
{Rates of convergence; the non-uniform case}
\setcounter{equation}{0}

In this section, we continue to build a theory  on the rate of
convergence of $EL$ for the non-uniform sample points which was
started by Hero,  Costa,  and Ma (2003). In Hero,  Costa,  and Ma
(2003) they studied the rate of convergence of $EL$ for the
non-uniform sample points  based on the very restrictive add-one
bound \eq{1.13} condition. Our starting point of this study is to
try to remove this condition in their argument and our work in
Section 2 is the result of this trial. In this section
 we build  a theory  on the rate of convergence of $EL$ for the
non-uniform sample points  which does not depend on the add-one
bound \eq{1.13} condition.

First, we work with a block density function; a probability
density function $\phi$ of the form
$$
\phi(x) = \sum_{i=1}^{m^d} \phi_i {\bold 1}_{Q_i} (x),
$$
where $\phi_i \ge 0$ is a constant and where $\{Q_i, 1 \leq i \leq
m^d\}$ is a partition  of $[0,1]^d$ into $m^d$ subboxes of edge
length $m^{-1}$, is  a block
density function of level $m$.\\

\noin {\bf Theorem 3.} Let $L$ be a strong Euclidean functional of
power $0<p<d$.  Then for i.i.d.\ sample points $\{ X_i : i
\geq 1 \}$ with common  block density function $\phi$ of level $m$
there exists a constant $C > 0$, independent of $\phi$, $n$, $m$,
such that \bea &&\left|
\frac{EL(\{X_1,\ldots,X_{n}\})}{n^{(d-p)/d}}-\alpha \int_{[0,1]^d}
\phi^{(d-p)/d}(x) \, dx\right| \nn\\
&\leq& \left\{ \ba{ll}
C(nm^{-d})^{-1/d}&\mb{for }0<p<d-1\\
C(nm^{-d})^{-(d-p)/d}\left(\log(nm^{-d})\vee 1\right)&\mb{for }p=d-1\neq 1\\
C(nm^{-d})^{-(d-p)/d}&\mb{for }d-1<p<d \\
C(nm^{-d})^{-(d-p)/d}&\mb{for }p=d-1=1,\ea \right. \lbl{3.1}\eea
where $\alpha$ is a constant  given by Theorem A. \ignore{A more
explicit form for the bound on the right hand side is
\begin{displaymath}
 O((nm^{-d})^{-1/d})= \left\{ \begin{array}{ll}
 \frac{C+C}{(nm^{-d})^{1/d}}\int f^{\frac{d-p-1}{d}}(x) dx
 (1+ o(1)),& d>2, \\
 \frac{C+C+C+\alpha}{(nm^{-d})^{1/d}}\int
 \phi^{\frac{d-p-1}{d}}(x) dx
 (1+ o(1)),& d=2. \end{array} \right.
\end{displaymath}
And in case $d-1 < p<d$, \be
\left|EL(X_1,\ldots,X_{n})/n^{(d-p)/d}-\alpha \int f^{(d-p)/d}(x) dx\right| \leq O((nm^{-d})^{-(d-p)/d}).
\lbl{3.2} \ee And in case $p=d-1$, \be
\left|EL(X_1,\ldots,X_{n})/n^{(d-p)/d}-\alpha \int f^{(d-p)/d}(x) dx\right| \leq
O(\log(\frac{n}{m^d})(nm^{-d})^{-(d-p)/d}). \lbl{3.2} \ee }

\noin {\bf Proof.} Hero,  Costa,  and Ma (2003) consider the case
$1\le p< d-1$ only. However, their argument also works for all the other three cases considered in this theorem
without any major changes.
For reader's convenience
here we reproduce the argument of Hero,  Costa,  and Ma (2003)
for the case $d-1<p<d$
and the case $p=d-1\neq 1$ with appropriate (minor) changes.

First, we handle the case  $d-1<p<d$. Let $X_1, \ldots, X_n$ be
$n$ i.i.d.\ samples with the common block density function
$\phi$ of level $m$ and let $L$ be a strong Euclidean functional.
Let $n_i$ be the number of samples $\{X_1,\ldots,X_n\}$ falling
into $Q_i$. Then, by \eq{1.5}, \eq{1.3}, Theorem 1 (i) \bea
EL(\{X_1,\ldots,X_n\},[0,1]^d,p)
&\leq&
\sum_{i=1}^{m^d}EL(\{X_1,\ldots,X_n\} \cap Q_i, Q_i, p) + C m^{d-p} \nn \\
&=&
m^{-p} \sum_{i=1}^{m^d} EE(L(\{U_1,\ldots,U_{n_i}\},[0,1]^d,p)|n_i) + Cm^{d-p}\nn\\
&\leq&
m^{-p} \sum_{i=1}^{m^d} E\left(\alpha n_i^{(d-p)/d} + C\right) + C m^{d-p}
\nn\\
&\le& \al m^{-p}
n^{(d-p)/d}\sum_{i=1}^{m^d}E\left(\frac{n_i}{n}\right)^{(d-p)/d} +
C m^{d-p}.\nn \eea By Jensen's inequality, for $0<\nu <1$ and $
p_i := \phi_i m^{-d}$ $$ E\left(\frac{n_i}{n}\right)^\nu \leq
p_i^\nu.$$ So, \bea
\frac{EL(\{X_1,\ldots,X_n\},[0,1]^d,p)}{n^{(d-p)/d}}
&\leq& \alpha \sum_{i=1}^{m^d} \phi_i^{(d-p)/d} m^{-d} + \frac{C}{(nm^{-d})^{(d-p)/d}}\nn \\
&=& \alpha \int_{[0,1]^d} \phi^{(d-p)/d} (x)\, dx
+\frac{C}{(nm^{-d})^{(d-p)/d}}. \nn \eea Similarly, by \eq{1.10}
\bea
& & EL^*(\{X_1,\ldots,X_n\},[0,1]^d,p) \nn \\
&\geq& m^{-p} \sum_{i=1}^{m^d} EE\left(L^*(\{U_1,\ldots,U_{n_i}\},[0,1]^d,p)|n_i\right) - C m^{d-p} \nn \\
&\geq& m^{-p} \sum_{i=1}^{m^d} EE\left(L(\{U_1,\ldots,U_{n_i}\},[0,1]^d,p)-C\Big|n_i\right) - Cm^{d-p} \nn \\
&\geq& m^{-p} \sum_{i=1}^{m^d} E\left(\alpha n_i^{(d-p)/d} - C-C\right)-C m^{d-p} \nn \\
&=& \al m^{-p} n^{(d-p)/d}
\sum_{i=1}^{m^d}E\left(\frac{n_i}{n}\right)^{(d-p)/d} -C
m^{d-p}.\nn \eea

Now, we claim that for $0<\nu<1$ and $p_i := \phi_i m^{-d}$ \be
E\left(\frac{n_i}{n}\right)^\nu \geq p_i^\nu - p_i^{\nu -
\frac{1}{2}}n^{-1/2}.\lbl{b} \ee If $g\in C^1(0,\infty)$ and if
$g$ is concave  over $x\ge 0$, monotone increasing over $x\ge 0$,
and $g(0) = 0$, then for any $x_0>0$,
$$ g(x) \geq g(x_0) - \frac{g(x_0)}{x_0}|x-x_0|. $$
Thus, with $g(x):= x^\nu$, $0<\nu<1$, $x:=n_i/n$, and $x_0:=p_i$,
we have
$$
\left(\frac{n_i}{n}\right)^\nu \geq p_i^\nu - p_i^{\nu -1}
\left|\frac{n_i}{n} - p_i\right|.$$ Take the expectation on both
sides. Since by Chebyshev's inequality
$$ E\left|\frac{n_i}{n} - p_i\right| \leq \left(E\left(\frac{n_i}{n} -
p_i\right)^2\right)^{1/2} \leq \frac{\sqrt{p_i}}{\sqrt{n}},$$
indeed we have \eq{b}.

So, by \eq{b}
\bea
&&\frac{EL^*(\{X_1,\ldots,X_n\},[0,1]^d,p)}{n^{(d-p)/d}}\nn\\
&\ge& \alpha \int_{[0,1]^d}\phi^{(d-p)/d}(x)\,dx-\frac{\alpha}{(nm^{-d})^{\frac{1}{2}}}
\int_{[0,1]^d} \phi^{\frac{1}{2}-\frac{p}{d}} (x)\,dx
-
\frac{C}{(nm^{-d})^{(d-p)/d}} . \nn \eea \ignore{ Since for $0<
\nu< 1$
$$
p_i^\nu - p_i^{\nu-\frac{1}{2}} /
\sqrt{n}
\le
E[(\frac{n_i}{n})^\nu] \leq p_i^\nu
$$
where $p_i:=\phi_i m^{-d}$ (see [3] for the detail), we have
\bea
E[L(\{X_1,\ldots,X_n\})/n^{(d-p)/d}]
&\leq& \alpha \sum_{j=1}^{m^d} \phi_i^{(d-p)/d} m^{-d} + \frac{C}{(nm^{-d})^{(d-p)/d}}\nn \\
&=& \alpha \int f^{(d-p)/d} (x) \,dx
+\frac{C}{(nm^{-d})^{(d-p)/d}}\nn \eea and \bea
& & E[L^*(X_1,\ldots,X_n)/n^{(d-p)/d}] \nn \\
&\geq& \alpha \sum_{j=1}^{m^d} \phi_i^{(d-p)/d} m^{-d}
-\frac{\alpha}{(nm^{-d})^{\frac{1}{2}}}
\sum_{j=1}^{m^d}\phi_i^{\frac{1}{2}-\frac{p}{d}} m^{-d} -\frac{C}{(nm^{-d})^{\frac{1}{d}}}\sum_{j=1}^{m^d}\phi_i^{\frac{d-1-p}{d}} m^{-d} \nn \\
& &  -
\frac{C}{(nm^{-d})^{(d-p)/d}}\nn \\
&=& \alpha \int_{[0,1]^d} f^{(d-p)/d}
(x)\,dx-\frac{\alpha}{(nm^{-d})^{\frac{1}{2}}} \int_{[0,1]^d}
f^{\frac{1}{2}-\frac{p}{d}} (x)\,dx
  - \frac{C}{(nm^{-d})^{(d-p)/d}}. \nn \eea
By the
definition of the Euclidean functional of power $p$
$$
E[L(X_1,\ldots,X_n)] \geq E[L^*(X_1,\ldots,X_n)] - Cn^{(d-p)/d} .
$$
} Therefore, by \eq{1.7} and \eq{1.10} we have
\bea
&&\left|\frac{EL(X_1,\ldots,X_{n})}{n^{(d-p)/d}}-\alpha \int_{[0,1]^d} \phi^{(d-p)/d}(x)\,dx\right|\nn\\
&\leq& \frac{\alpha}{(nm^{-d})^{1/2}}\int_{[0,1]^d} \phi^{\frac{1}{2}-\frac{p}{d}}(x)\,dx 
+  \frac{C}{(nm^{-d})^{(d-p)/d}}. \nn
\eea

Second, we handle the case $p=d-1\neq 1$ in a similar manner.
\bea
&&EL(\{X_1,\ldots,X_n\},[0,1]^d,p)\nn\\
&\leq& m^{-p} \sum_{i=1}^{m^d} E\left(\alpha n_i^{(d-p)/d} +C\left(\log n_i \vee 1\right) \right) + C m^{d-p}\nn\\
&\le& \al m^{-p}n^{(d-p)/d}\sum_{i=1}^{m^d}E\left(\frac{n_i}{n}\right)^{(d-p)/d} +Cm^{-p} \sum_{i=1}^{m^d}E\left(\log n_i \right) + C m^{d-p}. \nn
\eea
Since $n_i$ is highly concentrated around its mean $np_i$, for a large but fixed $C$ we have
$E \log n_i \le C \log n p_i$. So,
\bea
&&\frac{EL(\{X_1,\ldots,X_n\},[0,1]^d,p)}{n^{(d-p)/d}}\nn\\
 &\leq&
\alpha \sum_{i=1}^{m^d} \phi_i^{(d-p)/d} m^{-d} + \frac{C
\sum_{i=1}^{m^d} (\log \phi_i) m^{-d}}{(nm^{-d})^{(d-p)/d}} + \frac{C \log (nm^{-d})}{(nm^{-d})^{(d-p)/d}}
+ \frac{C}{(nm^{-d})^{(d-p)/d}}\nn \\
&\leq& \alpha \int_{[0,1]^d} \phi^{(d-p)/d} (x) \,dx+\frac{C
\int_{[0,1]^d} \log \phi (x) dx}{(nm^{-d})^{(d-p)/d}} + \frac{C
\log (nm^{-d})}{(nm^{-d})^{(d-p)/d}}
+\frac{C}{(nm^{-d})^{(d-p)/d}}. \nn \eea By the same reasoning,
\bea
&&\frac{EL^*(\{X_1,\ldots,X_n\},[0,1]^d,p)}{n^{(d-p)/d}}\nn\\
&\ge& \alpha \int_{[0,1]^d} \phi^{(d-p)/d}
(x)\,dx-\frac{\alpha}{(nm^{-d})^{1/2}}\int_{[0,1]^d}
\phi^{\frac{1}{2} -\frac{p}{d}}(x) \,dx-\frac{C \int_{[0,1]^d}
\log \phi (x) \,dx}{(nm^{-d})^{(d-p)/d}} \nn \\&-& \frac{C \log
(nm^{-d})}{(nm^{-d})^{(d-p)/d}} - \frac{C}{(nm^{-d})^{(d-p)/d}} .
\nn \eea Therefore, by \eq{1.7} and \eq{1.10} we have
\bea
& &\left|\frac{EL(X_1,\ldots,X_{n})}{n^{(d-p)/d}}-\alpha \int_{[0,1]^d} \phi^{(d-p)/d}(x)\,dx\right|  \nn \\
&\leq&\frac{C}{(nm^{-d})^{(d-p)/d}}\int_{[0,1]^d} \log \phi (x)\,dx +\frac{\alpha}{(nm^{-d})^{1/2}}\int_{[0,1]^d}\phi^{\frac{1}{2}-\frac{p}{d}}(x)\,dx\nn\\
&&\ \ \ \ \ +  \frac{C \log \frac{n}{m^{d}}}{(nm^{-d})^{(d-p)/d}} +
\frac{C}{(nm^{-d})^{(d-p)/d}}.\nn
\eea
\qed
\\

Second, we work with a probability density function $f$ in the
H$\ddot{o}$lder class $\sum (\beta,K,[0,1]^d)$, $\bt > 0$; a
probability density function $f:[0,1]^d \ra \mathR$ is in $\sum
(\beta,K,[0,1]^d)$ if for $x,y \in [0,1]^d$
$$
| f(y)-p_x^{\lfloor \beta \rfloor} (y)| \leq K|x-y|^\beta,
$$
where $\lfloor \beta \rfloor$ is the largest integer $l$ with
$l<\beta$ and where $p_x^{\lfloor \beta \rfloor}$ is the Taylor
expansion of $f$ at $x$ up to degree $\lfloor \beta \rfloor$.
\\

\noin{\bf Lemma 7.} If a probability density function $f$ is in
$\sum (\beta,K,[0,1]^d)$, $0<\beta \leq 1$, then for any partition
$\{Q_i, 1 \leq i \leq m^d\}$ of $[0,1]^d$ we define an
approximating block density function $\phi$ with level $m$ of the
form $ \phi(x) = \sum_{i=1}^{m^d} \phi_i {\bold 1}_{Q_i} (x) $ by
$\phi_i = m^d \int_{Q_i} f(x)dx.$ For this approximating block
density function $\phi$, \be \int_{[0,1]^d} \left|f(x) -
\phi(x)\right| \,dx \leq d^{\beta/2}K m^{-\beta}. \lbl{3.2}\ee

\noin {\bf Proof.}  By the mean value theorem there exists a point
$\eta_i \in Q_i$ such that
$$
\phi_i = m^d \int_{Q_i}f(x)\,dx = f(\eta_i).
$$
With this $\eta_i$, \bea \int_{[0,1]^d} \left|\phi(x)-f(x)\right|
\,dx &=&
\sum_{i=1}^{m^d} \int_{Q_i}  \left| f(\eta_i)- f(x) \right|\,dx\nn\\
&\le& \sum_{i=1}^{m^d} \int_{Q_i}  K\left| \eta_i-x
\right|^{\bt}\,dx\nn\\
&\le& d^{\bt/2}Km^{-\bt}.  \nn \eea
\qed\\

\noin{\bf Theorem 4.} Let $L$ be a strong Euclidean functional of
power $0<p<d$.  Then for i.i.d.\ sample points $\{ X_i : i
\geq 1 \}$ with common probability density function $f \in \sum
(\beta,K,[0,1]^d)$, $0<\beta \leq 1$, there exists a constant
$C:=C(d,\beta,p)
> 0$ such that
$$
\left|\frac{EL(X_1,\ldots,X_{n})}{n^{(d-p)/d}}-\alpha \int
f^{(d-p)/d}(x)\,dx\right| \leq \left\{ \ba{ll}
Cn^{-\frac{\beta (d-p)/d}{(\beta (d-p)/d+ 1) d }}&\mb{for }0<p<d-1\\
Cn^{-\frac{\beta}{d(\beta+d)}}\left(\log n\vee 1\right) &\mb{for }p=d-1\neq 1\\
Cn^{-\frac{\beta(d-p)/d}{\beta + d}}&\mb{for }d-1<p<d \\
Cn^{-\frac{\beta(d-p)/d}{\beta + d}} &\mb{for }p=d-1=1, \ea \right.
$$
where
$\alpha$ is a constant  given by Theorem A.

\ignore {One can easily modify the arguments in \cite{HCM03} to
prove theorem. In fact, one just need to replace certain estimates
in the proof of Proposition 2 and Proposition 3 in \cite{HCM03} by
the corresponding estimates in Theorem 3.
\qed \\

Last, we work with a probability density density function $f$ in
the  H$\ddot{o}$lder class $\sum (\beta,K,[0,1]^d)$, $\bt > 1$. To
do that we need to generalize the concept of the block density; a
probability density function $h$ of the form (with
$x:=(x_1,\ldots,x_d)$) \be h(x) = \sum_{i=1}^{m^d} \max \{
\sum_{\substack{(n_1,\ldots,n_d) \in
(\mathZ^{+})^d\\n_1+\cdots+n_d\le k}} a^{Q_i}_{n_{1}n_{2} \cdots
n_{d}} x_1^{n_1}\ldots x_d^{n_d},0 \} {\bold 1}_{Q_i} (x),
\lbl{3.A} \ee where $a^{Q_{i,0}}_{n_{1}n_{2} \cdots n_{d}} \ge 0$
and where $\{Q_i, 1 \leq i \leq m^d\}$ is a partition  of
$[0,1]^d$ into $m^d$ subboxes of edge length $m^{-1}$, is  a
$k$-th order block
density of $m$ level.\\

\ignore{
A
density $h(x)$ over $[0,1]^d$ is said to be the
\textit{$k$-th order density with level $m$} if for some set of
non-negative constants $\{a^{Q_{i}}_{n_{1}n_{2}
\cdots n_{d}}\}$ satisfying\\
$\sum_{i=1}^{m^d} \int_{Q_i} \max\{ \sum_{l=0}^{k}
\sum_{\substack{(n_1,\ldots,n_d)\\n_1+\cdots+n_d=l}}
a^{Q_{i}}_{n_{1}n_{2} \cdots n_{d}} x_1^{n_1}\ldots x_d^{n_d},0\}
dx = 1$,
$$ h(x) = \sum_{i=1}^{m^d} \max \{ \sum_{l=0}^{k}
\sum_{\substack{(n_1,\ldots,n_d)\\n_1+\cdots+n_d=l}}
a^{Q_{i}}_{n_{1}n_{2} \cdots n_{d}} x_1^{n_1}\ldots x_d^{n_d},0 \}
{\bold 1}_{Q_i} (x),$$ where $x = (x_1,\ldots,x_d)$ and $ \{ Q_i
\}_{i=1}^{m^d}$ is the uniform partition of the unit cube
$[0,1]^d$ defined previously.} }

\noin {\bf Proof.} As we did in  Lemma 7, for any partition
$\{Q_i, 1 \leq i \leq m^d\}$ of $[0,1]^d$ we define an
approximating
 block density $\phi$
with level $m$ of the form $ \phi(x) = \sum_{i=1}^{m^d} \phi_i
{\bold 1}_{Q_i} (x) $ by $\phi_i = m^d \int_{Q_i} f(x)\,dx$ for
which we have \eq{3.2}.

Here we consider the case $0<p<d-1$ only.  We leave all the other
cases to the reader as an exercise. Let $\{\tilde{X}_i\}_{i=1}^n$
be i.i.d.\ with the common probability density function  $\phi$.
Then,  \bea
& &\left|\frac{EL(\{X_1,\ldots,X_n\})}{n^{(d-p)/d}}-\alpha\int_{[0,1]^d} f^{(d-p)/d}(x)\,dx\right| \nn \\
&\leq& \left|\frac{EL(\{\tilde{X}_1,\ldots,\tilde{X}_n\})}{n^{(d-p)/d}} -\alpha \int_{[0,1]^d} \phi^{(d-p)/d} (x)\,dx\right|  \nn \\
&&+ \alpha \left|\int_{[0,1]^d} f^{(d-p)/d} (x)dx - \int_{[0,1]^d} \phi^{(d-p)/d} (x)\,dx\right| \nn \\
&&
+\frac{\left|EL(\{X_1,\ldots,X_n\})-EL(\{\tilde{X}_1,\ldots,\tilde{X}_n\})\right|}{n^{(d-p)/d}}
:= I + II + III. \nn \eea The term I can be handled by Theorem 3
and it is bounded by $C(nm^{-d})^{-1/d}$.

By \eq{3.2}, the term II is bounded by $C m^{-\beta(d-p)/d}$; \bea
II
&\leq& \alpha \int_{[0,1]^d} \left|f^{(d-p)/d}(x) - \phi^{(d-p)/d}(x)\right| \,dx \nn \\
&\leq& \alpha \int_{[0,1]^d} \left|f(x) - \phi(x)\right|^{(d-p)/d} \,dx \nn \\
&\leq& \alpha \left(\int_{[0,1]^d} \left|f(x) - \phi(x)\right| \,dx\right)^{(d-p)/d} \nn \\
&\leq& \alpha
\left(d^{\beta/2}K\right)^{(d-p)/d}m^{-\beta(d-p)/d}. \nn \eea

By standard coupling and by \eq{a}, the term III is bounded by
$Cm^{-\beta(d-p)/d}$;
\bea
 III
&=& \left|EL(\{X_1,\ldots,X_n\})-EL(\{\tilde{X}_1,\ldots,\tilde{X}_n\})\right|n^{-(d-p)/d} \nn \\
&\leq& C E\left|\{X_1,\ldots,X_n\}\triangle \{\tilde{X}_1,\ldots,\tilde{X}_n\}\right|^{(d-p)/d}n^{-(d-p)/d} \nn \\
&\leq& C \left(E\left|\{X_1,\ldots,X_n\}\triangle \{\tilde{X}_1,\ldots,\tilde{X}_n\}\right|\right)^{(d-p)/d} n^{-(d-p)/d}\nn \\
&\leq& C \left(E\sum_{i=1}^{n} {\bold 1}_{\{X_i \neq \tilde{X}_i\}}\right)^{(d-p)/d}n^{-(d-p)/d} \nn \\
&\leq& C\left(d^{\bt/2}K\right)^{(d-p)/d}m^{-\beta(d-p)/d}. \nn
\eea Therefore,
$$
\left|EL(X_1,\ldots,X_{n})/n^{(d-p)/d}-\alpha \int_{[0,1]^d}
f^{(d-p)/d}(x)\,dx\right| \leq \frac{C}{(nm^{-d})^{1/d}}+ C
m^{-\beta(d-p)/d}. $$

Now, we choose $m$ so that $(nm^{-d})^{-1/d}=m^{-\beta(d-p)/d}$,
i.e., we choose $m=n^{\frac{1}{\beta(d-p)+d}}$. Then, Theorem 4
follows from the above bounds for the terms I, II, and III.
\qed\\

It is an interesting problem to generalize Theorem 4 to cover more
general continuous probability density functions $f$. For example,
in the case $f \in \sum(\beta, K, [0,1]^d)$, $1< \beta \leq 2$, we
may work with a degree-one block density function; a probability
density function $\phi$ of the form
$$
\phi(x) = \sum_{i=1}^{m^d} \phi_i(x) {\bold 1}_{Q_i}(x),
$$
where $\phi_i(x) \geq 0$ is a degree-one function and where $\{
Q_i, 1 \leq i \leq m^d \}$ is a partition of $[0,1]^d$ into $m^d$
subboxes of edge length $m^{-1}$, is a degree-one block density
function of level $m$. For this degree-one block density function
using Theorem 4 we can get a rate of convergence similar to
Theorem 3. If we can get an approximation result using a
degree-one block density function $\phi$ similar to Lemma 7, we
can extend Theorem 4 to cover the case $1<\beta \leq 2$ (and
hopefully by iterating this argument) and the case
$2<\beta<\infty$. However, we face some technical problems in this
approximation procedure. We leave this problem to the interested
reader.

\vspace{1.0cm} \noin {\bf Acknowledgment.} The authors would like
to thank the referee for his thoughtful comments on an earlier
version of the paper. His valuable comments helped us to improve the
presentation of the paper.

\end{document}